\documentclass{amsart}
\usepackage{amsthm}
\usepackage{amsmath}
\begin{document}
\newcommand{\Q}{{\mathbb Q}}
\newcommand{\C}{{\mathbb C}}
\newcommand{\R}{{\mathbb R}}
\newcommand{\Z}{{\mathbb Z}}
\newcommand{\F}{{\mathbb F}}
\renewcommand{\wp}{{\mathfrak p}}
\renewcommand{\P}{{\mathbb P}}
\renewcommand{\O}{{\mathcal O}}
\newcommand{\Pic}{{\rm Pic\,}}
\newcommand{\Ext}{{\rm Ext}\,}
\newcommand{\rank}{{\rm rk}\,}
\newcommand{\sbull}{{\scriptstyle{\bullet}}}
\newcommand{\bX}{X_{\overline{k}}}
\newcommand{\ch}{\operatorname{CH}}
\newcommand{\tors}{\text{tors}}
\newcommand{\cris}{\text{cris}}
\newcommand{\alg}{\text{alg}}
\let\isom=\simeq
\let\rk=\rank
\let\tensor=\otimes

\newtheorem{theorem}{Theorem}[section]      
\newtheorem{lemma}[theorem]{Lemma}          %
\newtheorem{corollary}[theorem]{Corollary}  
\newtheorem{proposition}[theorem]{Proposition}
\newtheorem{scholium}[theorem]{Scholium}

\theoremstyle{definition}
\newtheorem{conj}[theorem]{Conjecture}
\newtheorem*{example}{Example}
\newtheorem{question}[theorem]{Question}

\theoremstyle{definition}
\newtheorem{remark}[theorem]{Remark}

\numberwithin{equation}{section}

\title{A $p$-adic proof of Hodge Symmetry for threefolds}
\author{Kirti Joshi}
\address{School of Mathematics, Tata Institute of Fundamental
Research, Homi Bhabha Road, Mumbai 400 005, INDIA and 
Math. department, University of Arizona, 617 N Santa Rita, Tucson
85721-0089, USA.}
\date{Preliminary Version: 9/9/2000}

\maketitle

\begin{quote}
\hfill Dedicated to the memory of P\=u.~La.~Deshpande.
\end{quote}


\section{Introduction}
Let $X/\C$ be a smooth projective variety. Then it is a consequence of
the Hodge decomposition theorem (see \cite{griffiths-harris}) that for
$p,q\geq 0$ we have:
\begin{equation}\label{hodge-symmetry}
h^{p,q}=\dim_{\C}H^q(X,\Omega^p_{X/\C})=h^{q,p}=\dim_{\C}
H^p(X,\Omega^q_{X/\C}).
\end{equation}
        The Hodge theorem also asserts that these two spaces are
complex conjugates and hence the equality of the two dimensions. In
this note we give a $p$-adic proof of \eqref{hodge-symmetry} when
$X/\C$ is a smooth projective threefold. Our approach is based on the
following observation: one first notes a purely $p$-adic assertion
that \eqref{hodge-symmetry} holds when Hodge numbers are replaced by
more delicate $p$-adic invariants introduced by Ekedahl
\cite{ekedahl84} called Hodge-Witt numbers.  These invariants take in
to account the torsion in the slope spectral sequence as well as the
slopes of Frobenius in the crystalline cohomology of the
variety. Hodge-Witt symmetry was proved by Ekedahl using his delicate
duality theorem (we note here that we do not use Ekedahl's duality in
the proof given below; we have replaced this by an elementary
assertion $h^{0,2}_W=h^{2,0}_W$) as well as the crystalline
Hard-Lefschetz theorem (when $X$ is a smooth projective threefold over
a perfect field of characteristic $p$).

        Once Hodge-Witt symmetry is proved one appeals to another
result of Ekedahl (see \cite{ekedahl-diagonal}) which guarantees the
equality of Hodge-Witt numbers and the Hodge numbers under suitable
circumstances. That these required conditions are met when $X$ is a
smooth projective variety over complex numbers is a simple consequence
of Deligne-Illusie criterion (see \cite{deligne87}) for degeneration
of Hodge-de Rham spectral sequence in characteristic $p$.

The restriction on dimension in Theorem~\ref{main} arises because
Hodge-Witt symmetry is not known.  Ekedahl gave a necessary and
sufficient condition for Hodge-Witt symmetry to hold (in any
dimension). This condition is given in terms of a certain equality of
domino numbers for the slope spectral sequence. While this condition
is not known to hold in dimension bigger than three, in
Remark~\ref{ekedahl-condition} and Scholium~\ref{ascholium} we give
one set of hypothesis under which these condition are satisfied.
Further Ekedahl's conditions also hold, for instance, if $X$ is a
smooth, projective Hodge-Witt variety over a perfect field of
characteristic $p>0$.  Hence we deduce that if a smooth projective
variety $X/\C$ has Hodge-Witt or ordinary reduction at an infinite set
of primes then Hodge symmetry holds in characteristic zero.

        We would like to thank Luc Illusie for correspondence and
comments; we would also like to thank Minhyong Kim for discussions on
Hodge symmetry and also to Jim Carlson whose lectures on mixed Hodge
theory at the University of Arizona revived our interest in algebraic
proofs of Hodge symmetry. We are grateful to F.~Oort for his interest
and his suggestions which have improved the readability of this
paper. We also thank Dinesh Thakur for comments.

\section{Hodge-Witt Symmetry}
        To keep this note brief, we will refer to \cite{illusie79b},
\cite{illusie83b} and \cite{ekedahl84}, \cite{ekedahl85},
\cite{ekedahl-diagonal} for notations and basic results.  In this
section $X/k$ is a smooth projective variety over a perfect field of
characteristic $p>0$. Let $H^j(X,W\Omega^i_X)$ be the Hodge-Witt
cohomology groups of $X$. Let $T^{i,j}$ be the dimension of the domino
associated to the differential $$H^j(X,W\Omega^i_X)\to
H^j(X,W\Omega^{i+1}_X).$$ Let $m^{i,j}$ be the slope numbers
associated to the slopes of Frobenius on the crystalline cohomology of
$X$. For the definition see \cite{illusie83a} or
\cite{ekedahl-diagonal}. Then the Hodge-Witt numbers of $X$, denoted
$h^{i,j}_W$ are defined to be
\begin{equation}\label{hodge-witt-numbers}
h^{i,j}_W=m^{i,j}+T^{i,j}-2T^{i-1,j+1}+T^{i-2,j+2}.
\end{equation}

Note that by \cite{illusie83b} $T^{i,j}$ is zero if the corresponding
differential of the slope spectral sequence is zero.

The following symmetry of slope numbers is a consequence of
\cite{deligne80} and \cite{katz74} and is due to Ekedahl (see
\cite{ekedahl-diagonal}) we give his proof here for completeness.
\begin{lemma}\label{slope-number-symmetry}
For any smooth projective variety $X/k$ over a finite field $k$ of
characteristic $p$ and for all $i,j$ we have
\begin{equation}
m^{i,j}=m^{j,i}.
\end{equation}
\end{lemma}
\begin{proof}
Let for any rational number $\lambda$ $h^n_{\cris,\lambda}$ be the
dimension (=multiplicity) of the slope $\lambda$ in
$H^n_{\cris}(X/W)$. Then by definition of $m^{i,j}$ we have 
\begin{equation}
m^{i,j}=\sum_{\lambda\in [i,i+1)}(i+1-\lambda)h^{i+j}_{\cris,\lambda}
        +\sum_{\lambda\in [i-1,i)}(\lambda-i+1)h^{i+j}_{\cris,\lambda}
\end{equation}
Thus to prove the equality of $m^{i,j}=m^{j,i}$ for all $i,j$ it
suffices to prove that if $\lambda$ occurs in $H^n_{\cris}(X/W)$ with
some multiplicity, then $n-\lambda$ occurs in this cohomology with the
same multiplicity. This is a consequence of Riemann hypothesis
\cite{katz74}, and here is how one derives it from the Riemann
Hypothesis. Let $\alpha_1,\ldots, \alpha_r$ be the eigen values of
Frobenius on $H^n_{\cris}(X/W)$. Then one knows by \cite{deligne80}
and \cite{katz74} that the characteristic polynomial of the Frobenius
on crystalline cohomology is the same as the characteristic polynomial
of Frobenius on $\ell$-adic \'etale cohomology for any $\ell\neq p$
(this is where we use the hypothesis that we are over a finite
field. In particular this characteristic polynomial has integer
coefficients and that the eigen values $\alpha_i$ are algebraic
integers with $\alpha_i\bar{\alpha_i}=q^n$, where $q$ is cardinality
of our finite ground field, and $q^n/\alpha_i$ is also an eigenvalue
of Frobenius on this cohomology. Normalise the $p$-adic valuation on
$W(k)[\alpha_1,\ldots,\alpha_r]$ so that $q$ has valuation $1$. Then
it follows from this that if $\lambda$ is the slope of frobenius then
$n-\lambda$ is also a slope of frobenius and it occurs with the same
multiplicity.
\end{proof}

We will also need the following elementary lemma due to Ekedahl (see
\cite{ekedahl-diagonal}).

\begin{lemma}\label{hodge-witt-symmetry}
Let $X/k$ be a smooth projective variety over a perfect field $k$ of
characteristic $p>0$. Then
\begin{equation}
h^{0,2}_W=h^{2,0}_W
\end{equation}
\end{lemma}
\begin{proof}
From the definition of $h^{0,2}_W$ and $h^{2,0}_W$ and
Lemma~\ref{slope-number-symmetry} it suffices to
prove that $T^{1,1}$ and $T^{2,0}$ are zero. The first is a consequence
of \cite{illusie83b} and the second is proved in \cite{illusie79b}. 
\end{proof}
\begin{remark}
One can give a similar proof of the fact that
$h^{0,1}_W=h^{1,0}_W=b_1/2$, where $\frac{1}{2}b_1$ is the dimension
of the Albanese variety of $X$ (see \cite{ekedahl-diagonal}).
\end{remark}

\begin{remark}
Ekedahl has shown that $h^{i,j}_W=h^{j,i}_W$ for any smooth projective
threefold $X/k$ and $i,j\geq 0$. The proof of the assertion
$h_W^{i,3-i}=h_W^{3-i,i}$ uses Ekedahl's duality theorem (see
\cite{ekedahl84}). Ekedahl's duality is not sufficient to prove the
corresponding assertion in higher dimensions. 
\end{remark}

\section{Hodge Symmetry for Threefolds}
\begin{theorem}\label{main}
Let $X/\C$ be a smooth projective variety over complex
numbers. Then we have $h^{i,j}=h^{j,i}$ for all $i+j\leq2$.
Moreover, if $X/\C$ is a smooth projective threefold then Hodge
Symmetry holds for $X$. 
\end{theorem}

\begin{proof}
Observe that $h^{i,3-i}=h^{3-i,i}$ is a trivial consequence of Serre
duality (see \cite{hartshorne-algebraic}). Hence the assertion for
threefolds follows from the first part of the assertion. So we have to
prove $h^{i,j}=h^{j,i}$ when $i+j\leq 2$. Observe that
$h^{0,1}=h^{1,0}$ follows by reduction to Picard variety of $X$, where
the assertion is trivial. So we are left with proving that
$h^{0,2}=h^{2,0}$.

We can assume, by a standard specialization argument, that $X$ is
defined over a number field. We choose a proper, regular model of $X$
over a suitable open set of ring of integers of the number field. By
further localisation on the base we may assume that model is smooth
over the base and that the relative de Rham cohomology is torsion free
and all the relative Hodge groups are torsion free as well. After
further localisation we can assume that the Hodge numbers of every
special fibre coincide with the Hodge numbers of the generic fibre and
that the Hodge to de Rham spectral sequence of the special fibre
degenerates at $E_1$. This can be done by the criterion of
degeneration of the Hodge de Rham spectral sequence due to
Deligne-Illusie (see \cite{deligne87}).

        Let $X_\wp$ be a special fibre chosen as above. Then
$H^*(X_\wp,W\Omega^\sbull_{X_\wp})$ is a Mazur-Ogus object in the
derived category of bounded complexes of modules over the
Cartier-Dieudonne-Raynaud algebra (see \cite{illusie83b}) and hence we
see that by \cite[Corollary~3.3.1, page 86]{ekedahl-diagonal} one sees
that $h^{i,j}_W=h^{i,j}$ and so by Lemma~\ref{hodge-witt-symmetry} we
are done.
\end{proof}

\begin{remark}\label{ekedahl-condition}
We note that of proof given above can be turned on its head: if $X$ is
a smooth projective variety over complex numbers then as one does know
that Hodge symmetry holds (\cite{griffiths-harris}), so for
sufficiently large primes, the reduction satisfies Hodge-Witt symmetry
as the Hodge numbers and the Hodge-Witt numbers coincide by
\cite{deligne87}, \cite{ekedahl-diagonal}. In particular one deduces
that Ekedahl's conditions on domino numbers (see
\cite[Proposition~3.2(ii), page 113]{ekedahl-diagonal}) are satisfied
for almost all but finite number of reductions, regardless of whether
or not these reductions are ordinary or Hodge-Witt. So it would seems
reasonable that at least under some reasonable hypothesis (liftablity
to $W_2$, torsion-free crystalline cohomology), Hodge-Witt symmetry
ought to hold in anay characteristic. We make this remark more precise
in the following.
\end{remark}

\begin{scholium}\label{ascholium}
  Let $X/k$ be a smooth projective variety over a perfect field of
  characteristic $p>0$. Assume that $X$ admits a smooth, projective
  lifting to $W(k)$ and $p>\dim(X)$ and the crystalline cohomology of
  $X$ is torsion free. Then for all $i,j$, the dimension, $T^{i,j}$ of
  the domino associated to the differential $d:H^j(X,W\Omega^i)\to
  H^j(X,W\Omega^{j+1})$ satisfies $T^{i,j}=T^{j-2,i+2}$.
\end{scholium}

\begin{remark}
  In a forthcoming work (see \cite{joshi2001}) we have investigated
  properties of Hodge-Witt numbers and slope numbers in detail. For
  instance we have shown that a weak version of the
  Bogomolov-Miyaoka-Yau inequality holds for all smooth projective
  surfaces of general type which lift to $W_2$ and have torsion free
  crystalline cohomology. We have also shown that the domino numbers
  of any smooth, projective variety which has torsion free crystalline
  cohomology are determined completely by their Hodge numbers and
  slope numbers. This was noted in \cite{ekedahl-diagonal} for Abelian
  varieties. Further we have also shown that for complete intersection
  in projective space, Ekedahl's conditions reduce to
  $T^{i,n-i}=T^{n-i-2,i+2}$ these are a consequnce of Ekedahl's duality
    for domino. These investigations will be reported in
    \cite{joshi2001}.
\end{remark}

\begin{remark}
Let $X/\C$ be  a smooth projective variety. Then  one expects that $X$
has Hodge-Witt or ordinary reduction  modulo an infinite set of primes
(see for instance \cite{joshi00a}). If  one believes the truth of these
assertions  then  by Ekedahl's  work  \cite{ekedahl-diagonal} and  the
above remarks, Hodge symmetry follows.
\end{remark}

The method of proof is more general and can be used to prove the
following:
\begin{theorem}\label{main2}
  Let $X/k$ be a smooth projective threefold over a perfect field $k$
  of characteristic $p>0$. Assume that the Hodge de Rham spectral
  sequence of $X$ degenerates and the crystalline cohomology of $X$ is
  torsion free. Then Hodge symmetry holds for Hodge numbers of $X/k$.
\end{theorem}

\begin{remark}
In general the conditions of Theorem~\ref{main2} are not very easy to
verify. But here is one application: assume that $X$ is Frobenius
split and that $p\geq 5$. Then by a result of Mehta (see
\cite{joshi00b}) we see that Hodge de Rham spectral sequence of $X$
degenerates at $E_1$, so if crystalline cohomology of $X$ is torsion
free then $X$ satisfies Hodge-Symmetry.
\end{remark}

\begin{remark}
        The assumption that the crystalline cohomology of $X$ is
torsion free is a necessary assumption in Theorem~\ref{main2}. If this
assumption is dropped then Hodge symmetry fails in positive
characteristic. Here is one example (for surfaces): this is taken from
\cite[Chapter II, Section 7.3]{illusie79b}. Let $X$ be a smooth
projective Enriques surface in characteristic two. Assume that $X$ is
singular, i.e., $H^1(X,\O_X)$ is one dimensional and frobenius is
bijective on this vector space; such surfaces exist only in
characteristic two. Then one has has a complete list of the Hodge
invariants of $X$. In particular in the present situation, there are
no global one forms on $X$ but as $H^1(X,O_X)\neq 0$, Picard scheme of
$X$ is not reduced (it is equal to $\mu_2$ so the Albanese variety is
zero and hence the second cristalline cohomology of X has torsion (it
is of type V-torsion in Illusie's classification of torsion). In fact
as the proof of the \cite[Proposition 7.3.5, page 656]{illusie79b}
shows, the cohomology of $X$ with coefficients in the sheaf of Witt
vectors is of finite type and so X is Hodge Witt.
\end{remark}


\begin{thebibliography}{10}

\bibitem{deligne80}
P.~Deligne.
\newblock La conjecture de {Weil} {II}.
\newblock {\em Publ. {Math}. {I.H.E.S}}, 52:137--252, 1980.

\bibitem{deligne87}
P.~Deligne and L.~Illusie.
\newblock Rel\'evements modulo $p^2$ et decomposition du complexe de de {R}ham.
\newblock {\em Invent. {M}ath.}, 89(2):247--270, 1987.

\bibitem{ekedahl84}
T.~Ekedahl.
\newblock On the multiplicative properties of the de {Rham}-{Witt} complex {I}.
\newblock {\em Ark. f\"ur {Mat}.}, 22:185--239, 1984.

\bibitem{ekedahl85}
T.~Ekedahl.
\newblock On the multiplicative properties of the de {Rham}-{Witt} complex
  {II}.
\newblock {\em Ark. f\"ur {Mat}.}, 23, 1985.

\bibitem{ekedahl-diagonal}
T.~Ekedahl.
\newblock {\em Diagonal complexes and {$F$}-guage structures}.
\newblock Travaux ex {C}ours. Hermann, Paris, 1986.

\bibitem{griffiths-harris}
P.~Griffiths and {J}. {H}arris.
\newblock {\em Principles of algebraic geometry}.
\newblock Wiley {C}lassics {L}ibrary. John {Wiley} and {Sons}, {Inc}., New
  {Y}ork, reprint of the 1978 edition edition, 1994.

\bibitem{hartshorne-algebraic}
R.~Hartshorne.
\newblock {\em Algebraic Geometry}.
\newblock Number~52 in Graduate {T}exts in {M}athematics. Springer-{V}erlag,
  New {Y}ork-{Heidelberg}, 1977.

\bibitem{illusie79b}
L.~Illusie.
\newblock Complexe de de {R}ham-{W}itt et cohomologie cristalline.
\newblock {\em Ann. Scient. Ecole Norm. Sup.}, 12:501--661, 1979.

\bibitem{illusie83a}
L.~Illusie.
\newblock {\em Algebraic {G}eometry {Tokyo/Kyoto}}, volume 1016 of {\em Lecture
  {Notes} in {M}athematics}, chapter Finiteness, duality and K\"unneth theorems
  in the cohomology of the de {Rham}-{Witt} complex, pages 20--72.
\newblock Springer-{Verlag}, 1983.

\bibitem{illusie83b}
Luc Illusie and {M}ichel {R}aynaud.
\newblock Les suites spectrales associ\'ees au complexe de de {R}ham-{W}itt.
\newblock {\em Inst. {H}autes {\'E}tudes {S}ci. {P}ubl. {M}ath.}, 57:73--212,
  1983.

\bibitem{joshi2001}
Kirti Joshi.
\newblock Crystalline aspects of geography of low dimensional varieties.

\bibitem{joshi00b}
Kirti Joshi.
\newblock Exotic torsion in {F}robenius split varieties.
\newblock Preprint, 2000b.

\bibitem{joshi00a}
Kirti Joshi and C.~S. Rajan.
\newblock On frobenius splitting varieties and ordinary varieties.
\newblock {\em Preprint.}, 2000.

\bibitem{katz74}
N.~Katz and W.~Messing.
\newblock Some consequences of the {R}iemann hypothesis for varieties over
  finite fields.
\newblock {\em Invent. {M}ath.}, 23:73--77, 1974.

\end{thebibliography}

\end{document}